\documentclass[11pt]{article}
\usepackage{definitions}
\usepackage{multirow}

\newcommand{\fperr}{\mathcal{E}}
\newcommand{\cupdlpx}{\texttt{cuPDLPx} }
\newcommand{\cupdlp}{\texttt{cuPDLP} }

\newcommand{\lc}{\ell_c}
\newcommand{\uc}{u_c}
\newcommand{\lv}{\ell_v}
\newcommand{\uv}{u_v}
\newcommand{\cR}{\mathcal{R}}
\newcommand{\cY}{\mathcal{Y}}

\newcommand{\cX}{\mathcal{X}}
\newcommand{\cS}{\mathcal{S}}
\newcommand{\R}{{\mathbb R}}

\title{cuPDLPx: A Further Enhanced GPU-Based First-Order Solver for Linear Programming}
\author{Haihao Lu\thanks{MIT, Sloan School of Management (haihao@mit.edu).} \and Zedong Peng\thanks{MIT, Sloan School of Management (zdpeng@mit.edu).} \and Jinwen Yang\thanks{University of Chicago, Department of Statistics (jinweny@uchicago.edu).}}
\date{}

\begin{document}
\maketitle

\begin{abstract}
We introduce \cupdlpx, a further enhanced GPU-based first-order solver for linear programming in C and CUDA. Building on the recently developed restarted Halpern PDHG for LP \cite{lu2024restarted}, \cupdlpx incorporates algorithmic techniques, including fixed-point-error-based restart criteria, a PID-controlled primal weight update, active-set stepsize boosting, and infeasibility detection. It also introduces a carefully engineered GPU C implementation that exploits modern CUDA features, including systematic kernel fusion, CUDA Graph execution, and an optimized sparse matrix-vector multiplication pipeline, together with integration of the PSLP presolver, all of which are major contributors to its practical performance. These improvements are carefully tailored for GPU architectures and deliver substantial computational gains. Compared with \texttt{cuPDLP.jl}, \cupdlpx \texttt{v0.4.0} achieves speedups of $7.1\times$--$8.4\times$ on MIPLIB 2017 LP relaxations and $9.9\times$--$11.5\times$ on Mittelmann's LP benchmark, and detects infeasibility or unboundedness on all 334 instances in the MIPLIB 2017 Infeas benchmark. The solver is publicly available at \url{https://github.com/MIT-Lu-Lab/cuPDLPx}.
\end{abstract}
\section{Introduction}
Linear programming (LP) is a central class of mathematical optimization problems due to its broad applicability, mathematical elegance, and computational efficiency~\cite{bertsimas1997introduction}. 
The two classic algorithmic frameworks for solving LPs are the simplex method~\cite{dantzig1948programming,dantzig1990origins} and interior-point methods (IPMs)~\cite{karmarkar1984new,wright1997primal}, which underpin the state of the art in high-accuracy CPU-based LP solvers. In broad terms, simplex methods traverse vertices of the feasible polyhedron and provide basic solutions, while IPMs follow a central path from the interior and offer polynomial-time guarantees. Despite decades of advances, both frameworks rely heavily on factorizations and triangular solves that are often challenging to scale to extremely large instances and to fully exploit modern massively parallel hardware.

{
Recently, first-order methods (FOMs) have emerged as an alternative for large-scale LPs~\cite{applegate2021practical,applegate2025pdlp}. In contrast to simplex and IPMs, FOMs avoid matrix factorizations and instead iterate using sparse matrix-vector multiplications (SpMV) and elementwise vector operations, making them particularly attractive when problem sizes are large. A notable example is PDLP~\cite{applegate2021practical,applegate2025pdlp}, which tailors restarted primal-dual hybrid gradient (PDHG) methods to LP and demonstrates strong practical performance via a set of enhancements.

The computational profile of FOMs also makes them a natural fit for modern parallel hardware, such as graphics processing units (GPUs). GPUs provide massive parallelism and high memory bandwidth, which align well with the dominant computations in PDHG-type methods. Motivated by this feature, recent work has developed GPU-accelerated LP solvers that adapt first-order frameworks to GPUs, including cuPDLP~\cite{lu2023cupdlp}, cuPDLP-C~\cite{lu2023cupdlpc}, HPR-LP~\cite{chen2024hpr}. These efforts highlight the practical potential of first-order LP solvers when tailored to modern hardware.
}

{ Building on these foundations, we present \cupdlpx, a GPU-accelerated solver for large-scale linear programming problems.
Compared with its predecessor \texttt{cuPDLP}, \cupdlpx incorporates recent algorithmic developments for first-order LP methods, including the restarted Halpern PDHG framework~\cite{lu2024restarted}. It also introduces a carefully engineered GPU C implementation that exploits modern CUDA features, including systematic kernel fusion, CUDA Graph execution, and an optimized sparse matrix-vector multiplication pipeline, all of which are major contributors to its practical performance. 
We summarize the primary contributions and positioning of \cupdlpx as follows:
\begin{itemize}
    \item \textbf{An efficient standalone GPU LP solver.}
    As a standalone software package, \cupdlpx provides an efficient GPU-based solver for large-scale linear programming problems. 
    The solver integrates algorithmic components motivated by recent theoretical and computational developments, including restarted Halpern PDHG, a constant step-size rule, a restart condition aligned with the analysis in~\cite{lu2024restarted}, and a primal-weight update based on a PID controller. 
    These components are combined with GPU-oriented implementation techniques, such as systematic kernel fusion and CUDA graphs, to improve end-to-end solver efficiency. 
    Most of the algorithmic ingredients are inherited from prior work~\cite{lu2024restarted}. 
    The contribution of \cupdlpx is therefore positioned at the solver-system level: it provides an efficient GPU solver implemented in C and evaluates the resulting end-to-end performance on standard LP benchmark sets. 
    For the first version of \cupdlpx, on the MIPLIB relaxations, \cupdlpx achieves a 1.3$\times$--1.9$\times$ speedup over \texttt{cuPDLP-C}; on Mittelmann's benchmark set, the speedups reach 1.9$\times$--3.1$\times$ across different tolerances and presolve settings.

    \item \textbf{An accessible template for GPU C implementation.}
    Beyond its performance as an LP solver, \cupdlpx is designed to provide a clean and easy-to-follow GPU C implementation of first-order optimization methods. 
    The implementation emphasizes lightweight data structures, modular solver components, and transparent GPU optimization techniques, in contrast to the previous C implementation \texttt{cuPDLP-C}~\cite{lu2023cupdlpc}. 
    This new design and open-sourced code lowers the barrier to developing efficient GPU-based first-order solvers and supports future extensions. 
    Its usefulness as an implementation reference is reflected in subsequent GPU-based optimization software and related implementations, such as~\cite{cuopt, li2026beyond, blin2026batched, li2026pdhcg}.
\end{itemize}

The solver is publicly available at \url{https://github.com/MIT-Lu-Lab/cuPDLPx}.}

{
\subsection{Related Work}
\label{sec:related_work}

{\textbf{ADMM-based solvers.}} A common first-order route is based on the alternating direction method of multipliers (ADMM). Representative solvers include OSQP~\cite{stellato2020osqp} and SCS~\cite{o2021operator}, which were developed primarily for quadratic and conic optimization and can be applied to LP through standard reformulations or embeddings. In these methods, each iteration typically involves solving a structured linear system, either through direct factorization of a regularized system matrix or through iterative linear solvers. Related LP-oriented approaches include ABIP~\cite{lin2021admm} and ABIP+~\cite{deng2022enhanced}, which combine ADMM-type iterations with homogeneous self-dual embedding and interior-point-inspired techniques.

\vspace{0.1cm}

\noindent{\textbf{Dual-ascent methods.}} Another line of work optimizes a smooth dual objective obtained from a regularized or perturbed LP.  ECLIPSE~\cite{basu2020eclipse} adds a small ridge regularization term to the primal LP, dualizes the complicating constraints, and applies first-order methods, including accelerated gradient methods, to the resulting smooth dual problem. More recently, DuaLip-GPU~\cite{dexter2026dualip} redesigns this ridge-regularized dual-ascent framework for GPU execution. 

\vspace{0.1cm}

\noindent{\textbf{PDHG-class solvers.}} A major recent direction for general-purpose large-scale LP is based on the primal-dual hybrid gradient method. PDLP~\cite{applegate2021practical,applegate2025pdlp} adapts restarted PDHG-type schemes to LP and incorporates practical enhancements such as diagonal preconditioning, adaptive restart, primal-weight adjustment, feasibility polishing, and infeasibility detection. Since each iteration primarily consists of sparse matrix-vector multiplications and simple projections, PDLP-style methods are naturally suited to parallel architectures. GPU implementations such as \texttt{cuPDLP}~\cite{lu2023cupdlp} and \texttt{cuPDLP-C}~\cite{lu2023cupdlpc} further demonstrate that GPU-oriented implementations of PDHG-class methods can help reduce time-to-solution on large LP instances. 

\vspace{0.1cm}

\noindent{\textbf{Restarted Halpern PDHG for LP.}}
Compared with the original PDLP framework~\cite{applegate2021practical,lu2023cupdlp,lu2023cupdlpc}, which employs restarted averaged PDHG, \cite{lu2024restarted} proposed restarted Halpern PDHG (rHPDHG) for solving LPs. The key algorithmic change is the replacement of the averaged restart mechanism with Halpern iteration and restarted based on subdifferential shrinkage, resulting in a more aggressive acceleration while maintaining numerical stability. From a theoretical perspective, rHPDHG directly controls the norm of a subdifferential residual rather than the primal--dual gap, leading to stronger guarantees, including accelerated two-stage convergence and faster infeasibility detection. Inspired by the later HPR-LP development~\cite{hprlp}, a reflected variant (called $r^2$HPDHG) together with its associated theoretical analysis was subsequently incorporated into an updated version of~\cite{lu2024restarted}. \cupdlpx is essentially based on the algorithmic and theoretical understanding of~\cite{lu2024restarted}. 

\vspace{0.1cm}

\noindent{\textbf{HPR-type methods.}} A related direction is based on Peaceman-Rachford-type operator splitting. HPR-LP~\cite{chen2024hpr} develops a GPU-based LP solver using a Halpern Peaceman-Rachford framework with semi-proximal terms. The method incorporates practical enhancements such as adaptive restarts and penalty-parameter updates, and reports strong empirical performance on standard LP benchmarks. 

These developments are closely related to \cupdlpx, and HPR-LP in particular provides an important reference point for GPU-based Halpern-type splitting methods for LP. Our work follows a complementary direction within the PDHG family.

Specifically, \cupdlpx is built upon the restarted Halpern PDHG framework~\cite{lu2024restarted}. In restarted Halpern PDHG, the primary iterates consist of the primal-dual variables, while auxiliary slack variables are recovered only when needed, such as for termination checks. In contrast, HPR-LP is derived from a Halpern Peaceman--Rachford splitting with semi-proximal terms, in which the corresponding slack variables are updated explicitly as part of each iteration.

Another important distinction concerns restart, which plays a central theoretical and practical role in LP. In \cupdlpx, restart is an integral component of the restarted Halpern PDHG framework and is motivated by achieving accelerated linear convergence~\cite{lu2024restarted}. In contrast, HPR-LP employs restart as a heuristic enhancement of the underlying HPR iteration, while its theoretical analysis establishes sublinear convergence~\cite{chen2024hpr}. Nevertheless, in practice, both approaches often exhibit linear convergence on most of the LP instances, and this behavior is a major contributor to their strong computational performance.

Beyond these differences, \cupdlpx was developed as a standalone C/CUDA solver with an emphasis on lightweight data structures, kernel fusion, CUDA Graph execution, and reduced memory movement. These implementation choices are central to the contribution of \cupdlpx, as they translate recent PDHG-class algorithmic developments into a practical C-based GPU LP solver.

As discussed in more detail below, \cupdlpx also incorporates several practical choices that are closely related to those used in HPR-LP, including the safeguard mechanism in the primal-weight update and vector normalization in the preconditioning step. 

}

\subsection{Paper organization}
Section \ref{sec:prelim} briefly introduces the form of LP to solve and the base algorithm vanilla PDHG. In Section \ref{sec:enhancement}, several algorithmic enhancements are proposed to accelerate convergence, followed by Section \ref{sec:impl}, introducing the implementation choices of \cupdlpx. The numerical comparisons between \cupdlpx and \cupdlp are presented in Section \ref{sec:numerical}.

\section{Preliminaries}\label{sec:prelim}

In this section, we introduce the LP form that \cupdlpx solves, followed by a discussion on vanilla PDHG for solving LPs.

\subsection{Notations}
For a symmetric positive semidefinite matrix $M \in \mathbb{R}^{n \times n}$, the {\( M \)-norm} of a vector \( x \in \mathbb{R}^n \) is defined as $\|x\|_M := \sqrt{x^\top M x}$. The orthogonal projector \( \operatorname{proj}_C(x) \) denotes the Euclidean projection of a point \( x \in \mathbb{R}^n \) onto a closed convex set \( C \subseteq \mathbb{R}^n\), i.e., $\operatorname{proj}_C(x) := \arg\min_{y \in C} \|x - y\|_2$. For a set \( C \subseteq \mathbb{R}^n \), we define \( -C := \{-x : x \in C\} \) as the {reflection of \( C \) about the origin}. The \((i,j)\)-th entry for a matrix \( M \in \mathbb{R}^{m \times n} \) is denoted by \( (M)_{ij} \). Given two vectors \( \ell \in (\mathbb{R} \cup \{-\infty\})^n \) and \( u \in (\mathbb{R} \cup \{\infty\})^n \) satisfying \( \ell \leq u \), define the mapping \( p: \mathbb{R}^n \times (\mathbb{R} \cup \{-\infty\})^n \times (\mathbb{R} \cup \{\infty\})^n \to \mathbb{R} \cup \{\infty\} \) given by $p(y; \ell, u) := u^\top y^+ - \ell^\top y^-$, {where $y^+ := \max\{y, 0\}$ and $y^- := \max\{-y, 0\}$
denote the componentwise positive and negative parts of $y$.}

\subsection{Linear programming}
\cupdlpx solves LP with the following primal-dual form:
\begin{equation}\label{eq:lp}
    \begin{aligned}[c]
    \min_{x\in \R^n}~~ &~ c^\top x \\
    \text{subject to:}~~ &~ \lc \leq Ax \leq \uc \\
    &~ \lv \leq x \leq \uv \ ,
    \end{aligned}
    \qquad\qquad
    \begin{aligned}[c]
    \max_{y\in \R^m,\, r\in \R^n} \quad & -p(-y; \lc, \uc) - p(-r; \lv, \uv) \\
    \text{subject to:}\quad & c - A^\top y = r \\
    & y\in\cY,~ r\in\cR \ ,
    \end{aligned}
\end{equation}

where $A \in \R^{m \times n}$, $c \in \R^{n}$, $\lc \in (\R \cup \{-\infty\})^m$, $\uc \in (\R \cup \{\infty\})^m$, $\lv \in (\R \cup \{-\infty\})^{n}$, $\uv \in (\R \cup \{\infty\})^{n}$, and the sets $\cY \subseteq \R^{m}$ and $\cR \subseteq \R^{n}$ are Cartesian products with their $i$th components given by
\begin{equation*}
\begin{aligned}
    &\cY_i := \begin{cases}
\{ 0 \} & (\lc)_i = -\infty, ~ (\uc)_i = \infty,  \\
\R^{-} & (\lc)_i = -\infty, ~ (\uc)_i \in \R ,\\
\R^{+} & (\lc)_i \in \R, ~ (\uc)_i = \infty, \\
\R & \text{otherwise}; 
\end{cases}
\quad \text{ and } \qquad
\end{aligned}
\begin{aligned}
    &\cR_i := \begin{cases}
\{ 0 \} & (\lv)_i = -\infty, ~ (\uv)_i = \infty,  \\
\R^{-} & (\lv)_i = -\infty, ~ (\uv)_i \in \R, \\
\R^{+} & (\lv)_i \in \R, ~ (\uv)_i = \infty, \\
\R & \text{otherwise};
\end{cases} \ .
\end{aligned}
\end{equation*}

For convenience, we define $\cX := \{x \in \R^n : \lv \leq x \leq \uv \}$ and  $\cS:=\{s\in\R^m: \lc\leq s \leq\uc\}$.
Equivalently, the primal-dual form of \eqref{eq:lp} is
\begin{equation}\label{eq:lp-pd}
\begin{aligned}
\min_{x\in\R^n}\max_{y\in\R^m}\quad & c^\top x-y^\top A x - p(y; -\uc, -\lc) \\
\text{subject to:}\quad & \lv\leq x\leq\uv,\quad y\in\cY.
\end{aligned}
\end{equation}

This LP form is used in CPU-based PDLP~\cite{applegate2025pdlp}.
\subsection{PDHG} 
\label{section:pdhg}
PDHG serves as the base routine of \cupdlpx to solve primal-dual problem \eqref{eq:lp-pd}. Specifically, the update rule of PDHG on \eqref{eq:lp-pd} is given as:
\begin{equation}\label{eq:pdhg}
    \begin{aligned}
        & x^{k+1} = \operatorname{proj}_{\cX}\pran{x^k-{\tau}(c-A^\top y^k)} \\
        & y^{k+1}=y^k-\sigma A(2x^{k+1}-x^k)-\sigma\operatorname{proj}_{-\cS}\pran{\sigma ^{-1}y^k-A(2x^{k+1}-x^k)}\ , 
    \end{aligned}
\end{equation}
where $\sigma > 0$ and $\tau > 0$ are dual and primal stepsizes respectively. \cupdlpx further reparameterizes as $\tau=\frac{\eta}{\omega}$ and $\sigma=\eta\omega$, where $\eta$ is called the stepsize and $\omega$ is the primal weight. And the canonical norm of PDHG is defined by $\|\cdot\|_P$, where $P:=P_{\eta,\omega}=\begin{bmatrix}
    \frac{\omega}{\eta}I & A^\top\\ A & \frac{1}{\eta\omega}I
\end{bmatrix} \succ 0$. For notational convenience, we define the primal-dual iterate at iteration $k$ as $z^k = \begin{bmatrix} x^k \\ y^k \end{bmatrix}$, and write $z^{k+1} = \operatorname{PDHG}(z^k)$ to denote a single PDHG update step applied to $z^k$ as defined in~\eqref{eq:pdhg}. 

{
Note that the reduced cost variable $r$ in \eqref{eq:lp} does not appear explicitly in the primal-dual form \eqref{eq:lp-pd} nor in the PDHG updates \eqref{eq:pdhg}. Hence, $r$ is not part of the iterates generated by PDHG and must be recovered from the primal-dual variables when evaluating the optimality and termination criteria~\cite{applegate2025pdlp}. 
Different recovery strategies are possible, and \cupdlpx recovers $r$ via the Moreau identity as $r^{k+1}={\tau}^{-1}(x^{k+1}-x^k)+(c-A^\top y^k)$.
}

\section{Algorithmic enhancements}\label{sec:enhancement}

In this section, we discuss a series of algorithmic further enhancements used in \cupdlpx to improve efficiency and convergence of the algorithm. These enhancements are primarily motivated by recent theoretical advances in restarted Halpern PDHG, as described in~\cite{lu2024restarted}.

At the high-level, \cupdlpx is based on the restarted Halpern PDHG algorithm (in contrast to the restarted average PDHG algorithm used in \cupdlp), and includes additional techniques such as the reflected updates, adaptive restarting, stepsize selection, active set stepsize boost, dynamic primal-dual weight, diagonal preconditioning and infeasibility detection. Each component is carefully designed and implemented to accelerate convergence and enhance stability for solving linear programming problems. We describe these enhancements below.

\subsection{Halpern scheme}
\label{section:halpern}
The Halpern scheme{~\cite{halpern1967fixed}}, along with its reflected variant, is an enhancement to the PDHG algorithm that improves convergence properties and accelerates solver performance. Originally developed to accelerate general operator splitting methods{~\cite{diakonikolas2020halpern,cai2022stochastic,tran2021halpern}}, the Halpern scheme has been successfully adapted to linear programming with strong theoretical guarantees and improved preliminary numerical performance~\cite{lu2024restarted}.
Halpern PDHG interpolates between the current PDHG iterate and the initial point, using a weighted average. Specifically, the update rule at iteration $k$ is:
\begin{equation*}
    z^{k+1}=\frac{k+1}{k+2}\operatorname{PDHG}(z^k)+\frac{1}{k+2}z^0\ .
\end{equation*}

As discussed in~\cite{lu2024restarted}, the Halpern and average schemes share structural similarities but differ in theoretical guarantees. Notably, the Halpern scheme supports reflection, which enables the use of larger stepsizes and often leads to faster empirical convergence, as described below.

\subsection{Reflection}
\label{section:reflection}
The reflection technique builds upon the Halpern scheme by applying a reflected version of the PDHG operator, $(1+\gamma)\operatorname{PDHG}-\gamma \operatorname{id}$ with reflection $\gamma \in [0,1]$, instead of the vanilla PDHG update. This is also known as over-relaxation in the literature. The resulting update rule of the Reflected Halpern PDHG update is
\begin{equation*}
    z^{k+1}=\frac{k+1}{k+2}\pran{(1+\gamma)\operatorname{PDHG}(z^k)-\gamma z^k}+\frac{1}{k+2}z^0\ .
\end{equation*}
This reflection mechanism effectively allows for longer steps than the vanilla Halpern update and has been shown, both theoretically and empirically, to improve convergence~\cite{lu2024restarted}. It is one of the key drivers of \cupdlpx’s superior performance compared to \cupdlp.

\subsection{Adaptive restart}
\label{section:restart}
Restarting is a key enhancement on first-order methods for attaining high-accuracy solutions for linear programming~\cite{applegate2023faster}. A restart is triggered when certain progress criteria are met, and the algorithm restarts from a new initial solution. In the Halpern scheme, the anchor is periodically reset to the current solution. This helps the algorithm stay focused on the neighborhood of the optimal solution, especially as the initial anchor becomes increasingly outdated during iterations. Such restart strategies have been shown to accelerate convergence rates for Halpern PDHG~\cite{lu2024restarted} in theory.

\cupdlpx adopts an adaptive restart strategy that evaluates potential restart conditions at each termination check. The core idea is to monitor a fixed-point error metric $\fperr(z)=\|z-\operatorname{PDHG}(z)\|_P$ at solution $z$, and trigger a restart when specific decay patterns are observed. This fixed-point error metric is motivated by the recent theoretical insight~\cite{lu2024restarted}, which shows that this fixed-point error is a natural metric for the restarted Halpern PDHG. It is different from the normalized duality gap used in CPU-based PDLP~\cite{applegate2021practical,applegate2025pdlp} and the KKT error used in cuPDLP~\cite{lu2023cupdlp}.

Three restart conditions are used:
\begin{itemize}
    \item \textbf{Sufficient decay.}
    \[
    \fperr(z^{n,k})\leq
    \beta_{\mathrm{sufficient}} \cdot \fperr(z^{n,0}).
    \]
    \item \textbf{Necessary decay and no local progress.}
    \[
    \begin{aligned}
    \fperr(z^{n,k})
    &\leq \beta_{\mathrm{necessary}} \cdot \fperr(z^{n,0})
    \quad\text{and} \quad
    \fperr(z^{n,k})
    &>\fperr(z^{n,k-1}).
    \end{aligned}
    \]
    \item \textbf{Artificial restart.}
    \[
    k\geq \beta_{\mathrm{artificial}} \cdot T,
    \]
    where $T$ is the total number of iterations.
\end{itemize}

This combined restart strategy balances aggressiveness and stability, contributing to cuPDLPx’s robust and accelerated convergence across a wide range of linear programming instances.

\subsection{Stepsize}
\label{section:stepsize}
In contrast to the adaptive stepsize heuristic used in \texttt{cuPDLP}, \cupdlpx initializes the stepsize at
\[
\eta^0=\frac{0.998}{\|A\|_2}
\]
, where the operator norm of the constraint matrix $\|A\|_2$ is approximated by power iteration. The stepsize remains fixed within each epoch. The rationale behind this change is twofold. First, the reflected Halpern PDHG inherently selects a larger effective stepsize, and the adaptive stepsize strategy from~\cite{lu2023cupdlp} can become overly aggressive in certain instances, leading to occasional instability. Second, this choice eliminates the need for the sequential stepsize search, which is not perfectly suited for parallel execution on GPUs. Thus, this modification improves both stability and parallel efficiency in \cupdlpx.

\textbf{Active-set stepsize boost.}
The initial stepsize $\eta^0$ is determined by the full constraint
matrix $A$. As the iterations progress, the active set may stabilize:
some variables remain clamped at their bounds and some constraints
remain inactive. Motivated by the finite identification properties of
restarted Halpern PDHG~\cite{lu2024restarted}, \cupdlpx uses an
estimate of this structure to form the reduced constraint matrix
$A_{B,F}$, where $B$ contains potentially binding constraints and $F$
contains variables that are not persistently clamped at a bound.
Since $\|A_{B,F}\|_2\le\|A\|_2$, its spectral norm may permit a
larger stepsize. This submatrix is used only for stepsize estimation and
all PDHG updates are applied to the full matrix $A$.
 
\textit{Active-set identification.}
At each termination check, the sets $B$ and $F$ are updated using
the underlying PDHG step~\eqref{eq:pdhg}, before
reflection and the Halpern update. This choice uses the active-set
information directly revealed by the projections, since the Halpern
update may change the active-set pattern of the projected point when
the anchor has a different active set.
To avoid reacting to transient active-set changes, an index is
excluded only if it passes the corresponding exclusion check below at
every termination check within a trailing window of $W$ iterations,
and is restored as soon as it fails the check.

For a current iterate $z=(x,y)$, let
$(\hat x,\hat y)=\operatorname{PDHG}(z)$.
Variable $i$ passes the exclusion check for $F$ if it is fixed,
$(\lv)_i=(\uv)_i$, or if either of the following reduced-cost checks
indicates that it is clamped at a finite bound:
\begin{equation}\label{eq:asb-var}
\begin{aligned}
c_i-(A^\top y)_i+\tau^{-1}(\hat x_i-x_i)
&>\epsilon_r\max\{1,|c_i|\}
&&\text{if }(\lv)_i>-\infty,\quad\text{or}\\
c_i-(A^\top y)_i+\tau^{-1}(\hat x_i-x_i)
&<-\epsilon_r\max\{1,|c_i|\}
&&\text{if }(\uv)_i<\infty.
\end{aligned}
\end{equation}
Here, $\epsilon_r>0$ is the tolerance for the reduced-cost sign check.

Constraint $j$ passes the exclusion check for $B$ if it is an inequality,
$(\lc)_j<(\uc)_j$, and its dual projection input lies strictly inside
$[-(\uc)_j,-(\lc)_j]$,
with a margin at each finite endpoint:
\begin{equation}\label{eq:asb-con}
\begin{aligned}
\sigma^{-1}y_j-[A(2\hat x-x)]_j+(\uc)_j
&>{\epsilon_c\max\{1,|(\uc)_j|\}}
&&\text{if }(\uc)_j<\infty,\quad\text{and}\\
-(\lc)_j-\sigma^{-1}y_j+[A(2\hat x-x)]_j
&>{\epsilon_c\max\{1,|(\lc)_j|\}}
&&\text{if }(\lc)_j>-\infty.
\end{aligned}
\end{equation}
Here, $\epsilon_c>0$ is the tolerance for identifying inactive constraints.
 
\textit{Stepsize boost.}
When boosting is enabled, \cupdlpx considers a stepsize increase at
adaptive restarts if at least one index has been excluded from $B$ or $F$.
The relative primal residual, relative dual residual, and relative
duality gap must also be below $\epsilon_{\mathrm{act}}$.
The estimate $\hat L\approx\|A_{B,F}\|_2$ is computed by masked
power iteration and reused until the retained sets change
sufficiently. The candidate stepsize is
\begin{equation}\label{eq:asb-target}
\eta^{+}
=\max\left\{\eta^0,\,
\min\left\{\frac{0.9}{\hat L},\,\eta_{\mathrm{ceil}}\right\}\right\}.
\end{equation}
Here, $0.9$ is the default safety factor, and
$\eta_{\mathrm{ceil}}$ is an upper limit on the boosted stepsize,
initially set to infinity and reduced after a rollback.
We set $\eta\gets\eta^+$ if
$\eta^+\ge\kappa\eta$, where $\kappa>1$ is the minimum increase ratio.
The stepsize remains fixed until the next boost or rollback.
 
\textit{Divergence protection.}
The stepsize boost relies on the identified active-set structure and
does not enforce the global PDHG stepsize condition $\eta\|A\|_2<1$.
If the active set changes, the boosted stepsize may become too large
and cause divergence.
During a boosted epoch, \cupdlpx uses the current PDHG operator
and the baseline metric $P^0:=P_{\eta^0,\omega}$ to evaluate the
fixed-point error for both adaptive restarting and rollback checks,
since $P_{\eta,\omega}$ may no longer be positive definite.
All compared errors, including that at the epoch anchor, are computed
using the stepsize $\eta$ and primal weight $\omega$ of the current epoch.
At each termination check during a boosted epoch, we trigger a
rollback if
\begin{equation}\label{eq:asb-divergence}
\left\|z^{n,k}-\operatorname{PDHG}(z^{n,k})\right\|_{P^0}
>
(1+\delta)
\left\|z^{n,0}-\operatorname{PDHG}(z^{n,0})\right\|_{P^0},
\end{equation}
where $\delta>0$ is the tolerance for the relative increase in the
fixed-point error.
Before each stepsize increase, the iterate and primal-weight
controller state are saved. A rollback first sets
$\eta_{\mathrm{ceil}}\gets\theta\eta$, where $\eta$ is the rejected
stepsize and $\theta\in(0,1)$.
It then restores the saved state, resets $\eta\gets\eta^0$, and starts
a new epoch. Boosting is disabled once the rollback limit is reached.

\subsection{Primal weight update}
\label{section:primal weight}
The primal weight $\omega$ is a crucial parameter in balancing progress between the primal and dual spaces. Our empirical experiences have shown that convergence behavior in PDLP is highly sensitive to the choice of this weight. The intuition for adjusting the primal weight $\omega^n$ is to set it such that the weighted distances to optimality in the primal and dual domains are in the same scale, i.e., 
$\|(x^{n,t} - x^*, 0)\|_{w^n} \approx \|(0, y^{n,t} - y^*)\|_{w^n}$. 
However, the optimal solution $(x^*, y^*)$ is not known during the iterations. 

To address this, \texttt{cuPDLP} estimates the primal and dual distances to optimality based on the observed movements in the previous epoch. The distance estimates are smoothed using exponential averaging and subsequently used to update the primal weight at each restart.
In \cupdlpx, we further enhance this strategy by modeling the update as a control problem and introducing a PID  controller to dynamically regulate the primal weight. 
We define the error as the gap between the primal and dual distances on a logarithmic scale:
\[
e^n = \log \left( \frac{\sqrt{w^n} \left\| x^{n,t} - x^{n,0} \right\|_2}{\frac{1}{\sqrt{w^n}} \left\| y^{n,t} - y^{n,0} \right\|_2}  \right)
\]

The primal weight is then updated according to:
\begin{equation}
    \label{eq:pid}
    \log w^{n+1} = \log w^n - [K_P \cdot e^n + K_I \cdot \sum_{i=1}^{n} \rho^{n-i} e^i + K_D \cdot (e^n - e^{n-1})]
\end{equation}
where $K_P$, $K_I$, and $K_D$ are the coefficients for the proportional, integral and derivative terms, respectively. {$\rho \in (0,1)$ is the decay factor for the accumulated errors. This logarithmic formulation has two advantages. First, applying the PID update in the logarithmic domain corresponds to a multiplicative update in the original scale of \(w\), yielding a scale-invariant controller. In other words, the control signal adjusts \(w\) by a relative factor rather than by an additive offset, making the update less sensitive to the absolute magnitude of \(w\). Second, since \(w\) is recovered from its logarithm, the positivity of the primal weight is automatically preserved.} The update is applied at each restart occurrence, and the initial primal weight is set to $1.0$.

{
To enhance numerical robustness, \cupdlpx also applies a safeguard rule for the primal-weight update, similar to HPR-LP~\cite{chen2024hpr}. 
When triggered, the safeguard overrides the PID update and resets the primal weight to the best historical value recorded so far (instead of resetting to $1$ as in HPR-LP).
Here, the best historical value is defined as the primal weight that has yielded the most balanced relative primal and dual residuals observed so far, as measured by
\[
\left|
\log_{10}
\frac{r_d^{\mathrm{rel}}}{r_p^{\mathrm{rel}}}
\right|,
\]
where \(r_p^{\mathrm{rel}}\) and \(r_d^{\mathrm{rel}}\) denote the relative primal and dual residuals, respectively.
This PID-controlled update, together with the safeguard, stabilizes the primal--dual balance across restarts and improves convergence speed and robustness.

\subsection{Diagonal preconditioner}
\label{section:diagonal preconditioner}
First-order methods such as PDHG are highly sensitive to the conditioning of the problem data. To mitigate the effects of ill-conditioning and improve numerical stability, \cupdlpx incorporates diagonal preconditioning by scaling the constraint matrix. Specifically, given the original constraint matrix $A$, a diagonal preconditioning step scales it to $\tilde A=D_1 A D_2$,
where $D_1$ and $D_2$ are positive diagonal matrices. This transformation yields a modified LP instance in which the matrix $\tilde A$ is better conditioned or ``well-balanced.'' The scaling also transforms the data as
: $\tilde{c} = D_2 c, \; (\tilde{l}_c, \tilde{u}_c) = D_1(l_c, u_c), \; (\tilde{l}_v,\tilde{u}_v) = D_2^{-1} (l_v,u_v)$.

\cupdlpx adopts the following diagonal preconditioning strategies, applied in the order listed below. The scaling factors are computed from the current matrix at each step, with unit factors for empty rows and columns:
\begin{enumerate}
    \item[(1)] Geometric mean scaling. Geometric mean scaling~\cite{tomlin1975scaling} scales each nonempty row and column by the reciprocal geometric mean of its smallest and largest nonzero entry magnitudes. At each step, $(D_1)_{jj}=\left(\min_{i:A_{ji}\ne0}|A_{ji}|\;\max_{i:A_{ji}\ne0}|A_{ji}|\right)^{-1/2}$ for each nonempty row $j$, and the column multipliers are computed analogously from the row-scaled matrix.
    \item[(2)] {Ruiz scaling.} Ruiz iterative scaling~\cite{ruiz2001scaling} scales both rows and columns of \( A \) to have similar norms. At each step, $(D_1)_{jj} = \|A_{j,\cdot}\|_\infty^{-1/2}, \; (D_2)_{ii} = \|A_{\cdot,i}\|_\infty^{-1/2}$ for $j = 1, \dots, m$ and $i = 1, \dots, n$ and repeated application drives the row and column norms toward 1.
    \item[(3)] {$\ell_1$ scaling.} The $\ell_1$ diagonal preconditioner is a special case of the general Pock-Chambolle diagonal scaling~ \cite{pock2011diagonal}. The diagonal matrices are defined as $(D_1)_{jj} = \|A_{j,\cdot}\|_{1}^{-1/2}, \; (D_2)_{ii} = \|A_{\cdot,i}\|_{1}^{-1/2}$ for $j = 1, \dots, m$ and $i = 1, \dots, n$.
    \item[(4)] {Vector normalization.} {The final step of diagonal preconditioning is to normalize the objective vector and constraint lower/upper bounds by their $\ell_2$ norms, respectively, similar to HPR-LP~\cite{hprlp}.} Specifically, $c$ is normalized by $\|c\|_2+1$, and $l_c, u_c, l_v, u_v$ are all divided by $\| [l_c, u_c] \|_{\text{eff}}+1$, where $\| \cdot \|_{\text{eff}}$ denotes a modified $\ell_2$ norm that only accounts for finite, non-duplicated boundary values.
    From the algorithmic perspective, it can be interpreted as an implicit normalization of the primal weight parameter.
\end{enumerate}

\subsection{Optimality termination criteria}
\label{section:optimality}
\cupdlpx terminates with an optimal solution when the relative KKT error is no greater than the termination tolerance $\epsilon>0$.
For a primal--dual point $(x,y)$ with recovered reduced cost $r$, this is equivalent to the following conditions on primal feasibility, dual feasibility, and the primal--dual gap:
\begingroup
\makeatletter
\@fleqnfalse
\makeatother
\begin{equation*}
\left\|Ax-\operatorname{proj}_{\cS}(Ax)\right\|_2
\leq\epsilon\left(1+\left\|(\lc,\uc)\right\|_2\right)
\tag*{\makebox[1pt][r]{(primal feasibility)}}
\end{equation*}
\begin{equation*}
\left\|c-A^\top y-r\right\|_2
\leq\epsilon(1+\|c\|_2)
\tag*{\makebox[1pt][r]{(dual feasibility)}}
\end{equation*}
\begin{equation*}
\left|c^\top x+p(-y;\lc,\uc)+p(-r;\lv,\uv)\right|
\leq\epsilon\left(1+\left|p(-y;\lc,\uc)+p(-r;\lv,\uv)\right|+|c^\top x|\right)
\tag{duality gap}
\end{equation*}
\endgroup
Infinite bound entries are treated as zero when computing $\|(\lc,\uc)\|_2$.
We check these conditions at the PDHG point before reflection and Halpern anchoring, using the original, unscaled problem data.
The reduced cost $r$ is recovered as described in Section~\ref{section:pdhg} and satisfies $r\in\cR$.

\subsection{Infeasibility detection}
\label{section:infeasibility}
 
In addition to identifying optimal solutions, \cupdlpx detects
primal and dual infeasibility through candidate rays extracted
from the restarted Halpern PDHG iterates.
For infeasible LPs, the differences of consecutive PDHG iterates and the
normalized iterates converge to the infimal displacement vector, whose nonzero
primal or dual component provides an infeasibility certificate
\cite{applegate2024infeasibility}.
For restarted Halpern PDHG, Lu and Yang~\cite{lu2024restarted} establish sublinear
convergence of the fixed-point difference and normalized displacement from
the restart anchor within each epoch, together with linear convergence
of the fixed-point difference at restart points.
 
\textbf{Candidate rays.}
Motivated by these results, at each termination check \cupdlpx examines
three candidate rays constructed from the available PDHG iterates:
\begin{itemize}
\item \textit{Fixed-point difference:}
$\operatorname{PDHG}(z^{n,k})-z^{n,k}$.
\item \textit{Current iterate:}
$\operatorname{PDHG}(z^{n,k})$.
\item \textit{Displacement from the epoch anchor:}
$\operatorname{PDHG}(z^{n,k})-z^{n,0}$.
\end{itemize}
Here, $z^{n,0}$ is the anchor of the current epoch.
These candidates capture complementary information: the fixed-point
difference benefits from the accelerated convergence associated with
restarting, the displacement from the anchor captures progress within
the current epoch, and the current iterate retains information accumulated
across previous epochs.

\textbf{Certificate test.}
For a candidate ray $(x,y)$, Farkas' lemma gives the following
conditions for certifying infeasibility.

\textit{Dual infeasibility.}
A primal ray $x$ certifies dual infeasibility if
\begin{equation}\label{eq:dual-infeas-certificate}
x\in\cX_\infty,\qquad Ax\in\cS_\infty,\qquad c^\top x<0,
\end{equation}
where $\cX_\infty$ and $\cS_\infty$ denote the recession cones
of $\cX$ and $\cS$, respectively.

\textit{Primal infeasibility.}
A dual ray $(y,r)$ certifies primal infeasibility if
\begin{equation}\label{eq:primal-infeas-certificate}
y\in\cY,\qquad r=-A^\top y\in\cR,\qquad
-p(-y;\lc,\uc)-p(-r;\lv,\uv)>0.
\end{equation}

To check these conditions numerically, \cupdlpx first projects the
primal and dual components of each candidate ray $(x,y)$ onto the
corresponding sign constraints:
$\bar x=\operatorname{proj}_{\cX_\infty}(x)$ and
$\bar y=\operatorname{proj}_{\cY}(y)$.
We also set $\bar r=\operatorname{proj}_{\cR}(-A^\top\bar y)$.

\begin{samepage}
\textit{Primal infeasibility criteria.}
\cupdlpx detects primal infeasibility if any one of the three candidate rays
satisfies both conditions:
\begin{subequations}\label{eq:primal-infeas-criterion}
\begin{equation}\label{eq:primal-infeas-violation}
\frac{\left\|A^\top\bar y+\bar r\right\|_\infty}
{\left\||A|^\top|\bar y|\right\|_\infty}
\le\epsilon_{\mathrm{infeas}},
\end{equation}
\begin{equation}\label{eq:primal-infeas-margin}
\frac{-p(-\bar y;\lc,\uc)-p(-\bar r;\lv,\uv)}
{p(-\bar y;-|\lc|,|\uc|)+p(-\bar r;-|\lv|,|\uv|)}
\ge\delta_{\mathrm{infeas}},
\end{equation}
\end{subequations}
where $|\cdot|$ is taken entrywise, and terms involving
infinite bounds are omitted.
\end{samepage}

\begin{samepage}
\textit{Dual infeasibility criteria.}
\cupdlpx detects dual infeasibility if any one of the three candidate rays
satisfies both conditions:
\begin{subequations}\label{eq:dual-infeas-criterion}
\begin{equation}\label{eq:dual-infeas-violation}
\frac{\left\|A\bar x-\operatorname{proj}_{\cS_\infty}(A\bar x)\right\|_\infty}
{\left\||A|\,|\bar x|\right\|_\infty}
\le\epsilon_{\mathrm{infeas}},
\end{equation}
\begin{equation}\label{eq:dual-infeas-margin}
\frac{-c^\top\bar x}{|c|^\top|\bar x|}
\ge\delta_{\mathrm{infeas}}.
\end{equation}
\end{subequations}
\end{samepage}

Conditions \eqref{eq:primal-infeas-violation} and
\eqref{eq:dual-infeas-violation} bound the relative violations of the
certificate constraints by $\epsilon_{\mathrm{infeas}}$, while
\eqref{eq:primal-infeas-margin} and \eqref{eq:dual-infeas-margin}
require a normalized objective margin of at least $\delta_{\mathrm{infeas}}$.
This separation of the two checks contrasts with the criteria used by
PDLP~\cite{applegate2021practical}, which bound the ratio of the
constraint violation to the objective term.
PDLP's criteria are invariant to positive scalar scaling of a
candidate ray, but not to such scaling of the objective coefficients or bounds.
For example, for a fixed candidate primal ray, multiplying $c$ by $t>0$
reduces the ratio used in PDLP's dual infeasibility check by a factor of $t$,
without changing the constraint violation.
Similarly, for a fixed candidate dual ray, multiplying all finite constraint
and variable bounds by $t$ reduces the ratio used in PDLP's primal
infeasibility check by the same factor.
However, our primal and dual infeasibility criteria,
\eqref{eq:primal-infeas-criterion} and \eqref{eq:dual-infeas-criterion},
remain unchanged under these scalar scalings.
The objective margins in \eqref{eq:primal-infeas-margin} and
\eqref{eq:dual-infeas-margin} are also invariant under positive diagonal
scaling for corresponding candidate rays. However, the relative constraint violations in
\eqref{eq:primal-infeas-violation} and \eqref{eq:dual-infeas-violation}
may depend on relative column and row scaling, respectively,
because of their global $\ell_\infty$ normalization.
\cupdlpx therefore evaluates these two constraint violations on the
diagonally preconditioned problem of
Section~\ref{section:diagonal preconditioner}, where the rows
and columns are equilibrated.

\subsection{Summary}
An overview of the restarted Halpern PDHG algorithm and the proposed enhancements implemented in \cupdlpx is summarized in Alg. \ref{alg:overview}.

\begin{algorithm}[htb]
\caption{Overview of the core algorithm implemented in \cupdlpx.}
\label{alg:overview}
\KwIn{Linear program \eqref{eq:lp}, initial point $(x^{0,0},y^{0,0})$, termination tolerance $\epsilon$.}
Diagonal preconditioning: geometric mean scaling, Ruiz scaling, $\ell_1$~scaling, and vector normalization \tcp*{Section~\ref{section:diagonal preconditioner}}
Initialize: $\eta^0 \gets 0.998 / \| A \|_2$, $\eta \gets \eta^0$, $w^0 \gets 1$ \tcp*{Sections~\ref{section:stepsize} and~\ref{section:primal weight}}
\For{$n \gets 0, 1, \dots$}{
    $k \gets 0$\;
    $\tau \gets \eta/w^n$, $\sigma \gets \eta w^n$\;
    \Repeat{restart condition holds \tcp*{Section~\ref{section:restart}}}{
        $\hat x^{n,k+1} \gets \operatorname{proj}_{X}(x^{n,k}-\tau (c-A^\top y^{n,k}))$ \tcp*{Section~\ref{section:pdhg}}
        $\hat y^{n,k+1} \gets y^{n,k}-\sigma A(2 \hat x^{n,k+1}-x^{n,k}) - \sigma \operatorname{proj}_{[-u_c,-l_c]}\left(\sigma ^{-1}y^{n,k}-A(2 \hat x^{n,k+1}-x^{n,k})\right)$\;
        \If{$\mathrm{KKT}(\hat x^{n,k+1},\hat y^{n,k+1})\le\epsilon$ \textbf{or} infeasibility detected \tcp*{Sections~\ref{section:optimality} and~\ref{section:infeasibility}}}{
            \Return $(\hat x^{n,k+1},\hat y^{n,k+1})$ or infeasibility certificate\;
        }
        $x^{n,k+1}=\frac{k+1}{k+2}\pran{(1+\gamma)\hat x^{n,k+1}-\gamma x^{n,k}}+\frac{1}{k+2} x^{n,0}$ \tcp*{Sections~\ref{section:halpern} and~\ref{section:reflection}}
        $y^{n,k+1}=\frac{k+1}{k+2}\pran{(1+\gamma)\hat y^{n,k+1}-\gamma y^{n,k}}+\frac{1}{k+2} y^{n,0}$\;
        $k \gets k+1$\;
    }
    $(x^{n+1,0}, y^{n+1,0}) \gets (\hat x^{n,k}, \hat y^{n,k})$\;
    Update primal weight $w^{n+1}$ \tcp*{Section~\ref{section:primal weight}}
    Apply the active-set stepsize boost if eligible \tcp*{Section~\ref{section:stepsize}}
}
\end{algorithm}

\section{Implementation choices of \cupdlpx}\label{sec:impl}

This section summarizes the main GPU-oriented implementation choices in \cupdlpx, in particular stable and efficient SpMV, systematic kernel fusion, CUDA Graph execution, and presolve integration. These implementation choices constitute an important part of the software contribution and have a significant impact on the solver's practical performance.

\subsection{Sparse Matrix--Vector Multiplication}

SpMV is the dominant kernel in \cupdlpx. 
To avoid the overhead of transposed SpMV calls, \cupdlpx explicitly stores both $A$ and $A^\top$ in compressed sparse row (CSR) format. 
This increases matrix storage but allows both $Ax$ and $A^\top y$ to be computed using non-transposed SpMV kernels, which are generally more efficient in \texttt{cuSPARSE}.

For CUDA versions before 13.3, \cupdlpx uses the \texttt{cuSPARSE\_SpMV} API with \texttt{CUSPARSE\_SPMV\_CSR\_ALG2}. 
This algorithm provides deterministic results across runs, which is useful for reproducibility. 
We also call \texttt{cusparseSpMV\_preprocess()} so that the one-time preprocessing cost can be amortized over the many SpMV calls made during the solve. 
Starting from CUDA~13.3, \cupdlpx additionally supports \texttt{cusparseSpMVOp()} as an alternative backend. 
On our LP workloads, this backend is often faster than the deterministic \texttt{cusparseSpMV()} configuration, at the cost of a larger temporary buffer.

\subsection{Kernel Fusion}

Besides SpMV, most operations in \cupdlpx are bandwidth-bound vector primitives, including affine updates, scaling, clipping, and projections. 
To reduce global-memory traffic, \cupdlpx fuses consecutive elementwise operations whenever possible. 
In iterations without termination or restart checks, the intermediate variables $\hat x^{n,k+1}$ and $\hat y^{n,k+1}$ are not written to global memory, and the dual slack variable $r$ is not formed. 
Thus, each iteration consists of two SpMV calls and two fused update kernels, one for the primal update and one for the dual update.

When termination or restart checks are performed, the required intermediate quantities and slack variables are computed within the same fused update kernels.  
All work buffers are preallocated during initialization and reused throughout the solve, avoiding dynamic allocation in the main iteration loop.

\subsection{CUDA Graphs}

Although kernel fusion reduces memory traffic, first-order methods still invoke many short GPU kernels. 
The resulting kernel-launch overhead can become a nontrivial part of total runtime. 
To reduce this overhead, \cupdlpx uses CUDA Graphs: the repeated GPU workflow is captured once, instantiated, and then replayed many times. 
This reduces CPU-side launch overhead and may also enable the CUDA compiler to optimize the captured workflow.

Rather than capturing only a single iteration, \cupdlpx captures a multi-iteration window aligned with the termination and restart check frequency. 
This amortizes graph-instantiation overhead and increases the benefit of graph replay.

\subsection{Presolve}

\cupdlpx integrates the PSLP presolver~\cite{Cederberg26} before the main PDHG iterations. PSLP applies standard LP reductions, such as removing redundant constraints, tightening bounds, and fixing variables when possible, to produce a smaller equivalent problem. It consists of three stages: initialization of presolve data structures, presolve reductions, and postsolve recovery of a primal--dual solution for the original LP.  In addition to reducing problem size, presolve can improve conditioning and stabilize the subsequent first-order iterations.
}

\section{Numerical experiments}\label{sec:numerical}
In this section, we evaluate the numerical performance of \cupdlpx in solving LPs to optimality and detecting infeasibility. We first describe the experimental setup and compare \cupdlpx \texttt{v0.1.0} with other GPU-based solvers on MIPLIB 2017 LP relaxations and Mittelmann's LP benchmark. Section~\ref{sec:version-comparison} compares successive releases of \cupdlpx, and Section~\ref{sec:infeasibility-experiments} evaluates infeasibility detection on two additional benchmark sets.

\textbf{Benchmark datasets.}
For the optimality experiments, we use two widely used LP benchmark datasets: MIP Relaxations, which consists of 383 instances derived from root-node LP relaxations of MIPLIB 2017~\cite{gleixner2021miplib}, and Mittelmann’s LP benchmark dataset~\cite{mittelmannbenchmark}, which contains 49 LP instances. In particular, the 383 instances in the MIP relaxations dataset are selected based on the same filtering criteria outlined in~\cite{lu2023cupdlp}.
For the infeasibility experiments, we use the 334 instances of MIPLIB 2017 Infeas~\cite{lu2026hpdlp} and the 29 infeasible Netlib instances~\cite{netlib}. MIPLIB 2017 Infeas comprises 229 primal-infeasible LPs and 105 feasible but unbounded LPs derived from MIPLIB 2017.
Table~\ref{tab:benchmark-summary} summarizes the instance counts and problem size ranges for these four benchmark sets.
Based on the number of non-zeros of the constraint matrices, we further categorize the 383 MIP Relaxations instances into three groups, as summarized in Table~\ref{tab:miplib-size}.
\begin{table}[ht!]
\centering
\begin{tabular}{lcccc}
\toprule
\textbf{Benchmark set} & \textbf{Instances} & \textbf{Rows} & \textbf{Columns} & \textbf{Nonzeros} \\
\midrule
MIP Relaxations & 383 & 36--19.9M & 150--20.7M & 100.9k--257.0M \\
Mittelmann & 49 & 4.3k--30.0M & 6.6k--38.9M & 95.0k--93.0M \\
MIPLIB 2017 Infeas & 334 & 27--19.9M & 297--20.7M & 100.9k--58.6M \\
Netlib (infeasible) & 29 & 8--3.8k & 4--10.7k & 16--97.2k \\
\bottomrule
\end{tabular}
\caption{Number of instances and ranges of problem sizes for each benchmark set.\label{tab:benchmark-summary}}
\end{table}

\begin{table}[ht!]
\centering
\begin{tabular}{cccc}
\hline
                    & \textbf{Small}                   & \textbf{Medium}                     & \textbf{Large}                   \\ \hline
\textbf{Number of nonzeros}  & 100K -  1M & 1M - 10M & \textgreater 10M \\
\textbf{Number of instances} & 269                     & 94                         & 20                      \\ \hline
\end{tabular}
\caption{Scale and number of instances in {MIP Relaxations}.\label{tab:miplib-size}}
\end{table}

\textbf{Software.} 
We implement \cupdlpx in C and directly leverage low-level CUDA libraries for improved performance. In the optimality comparisons in Tables~\ref{tab:miplib-no-presolve}--\ref{tab:speedup-mittelmann}, we compare \cupdlpx \texttt{v0.1.0} against \texttt{cuPDLP.jl}\footnote{\texttt{cuPDLP.jl}: \url{https://github.com/jinwen-yang/cuPDLP.jl}, accessed September 2025}, \texttt{HPR-LP v0.1.0}\footnote{\texttt{HPR-LP v0.1.0}: \url{https://github.com/PolyU-IOR/HPR-LP}, accessed August 2025}, \texttt{Gurobi 12.0.4}\footnote{\texttt{Gurobi 12.0.4}: \url{https://www.gurobi.com}, accessed September 2025}, and \texttt{cuPDLP-C v0.4.1}\footnote{\texttt{cuPDLP-C v0.4.1}: \url{https://github.com/COPT-Public/cuPDLP-C}, accessed September 2025}. Both Julia implementations, \texttt{cuPDLP.jl} and \texttt{HPR-LP}, are timed after precompilation to exclude Julia compilation overhead.
For the Gurobi-PDLP results in these optimality comparisons, we set \texttt{Method=6} to activate PDLP and \texttt{GURO\_PAR\_PDHGGPU=1} to enable GPU acceleration.

The infeasibility experiments use \cupdlpx\ \texttt{v0.4.0}.
The baselines are \texttt{Gurobi}~\cite{gurobi} with its primal simplex, dual simplex, barrier, and homogeneous barrier methods, together with
\texttt{HiGHS}\footnote{\texttt{HiGHS}: \url{https://github.com/ERGO-Code/HiGHS}, accessed October 2026}~\cite{huangfu2018parallelizing},
\texttt{SCS}\footnote{\texttt{SCS}: \url{https://www.cvxgrp.org/scs/}, accessed October 2026}~\cite{o2021operator},
\texttt{OR-Tools PDLP}\footnote{\texttt{OR-Tools PDLP}: \url{https://github.com/google/or-tools}, accessed October 2026}~\cite{applegate2025pdlp}, and
\texttt{H-PDLP}\footnote{\texttt{H-PDLP}: \url{https://github.com/MIT-Lu-Lab/H-PDLP}, accessed October 2026}~\cite{lu2026hpdlp}.
Unless otherwise specified, all solvers use their default parameter settings in these infeasibility experiments.

\textbf{Computing environment.} 
We use an NVIDIA H100-SXM-80GB GPU with CUDA 12.4. The
experiments of \texttt{cuPDLP.jl} and \texttt{HPR-LP} are performed in Julia 1.11.5.

\textbf{Initialization.} 
\cupdlpx uses all-zero vectors as the initial starting points.

\textbf{Optimality tolerances.} \cupdlpx uses the termination criteria in Section~\ref{section:optimality}. In the optimality experiments, we set $\epsilon=10^{-4}$ for moderate accuracy and $\epsilon=10^{-8}$ for high accuracy.
For Gurobi-PDLP, both absolute and relative feasibility tolerances are checked, and feasibility is considered satisfied once either tolerance is met.
To align termination criteria across solvers, we fix the absolute feasibility tolerance to \texttt{GURO\_PAR\_PDHGABSTOL=1e-9}, while \texttt{GURO\_PAR\_PDHGCONVTOL} and \texttt{GURO\_PAR\_PDHGRELTOL} serve as the effective convergence criteria.

\textbf{Infeasibility termination criteria.}
\cupdlpx uses the certificate tests described in Section~\ref{section:infeasibility}, with $\epsilon_{\mathrm{infeas}}=10^{-8}$ and $\delta_{\mathrm{infeas}}=10^{-8}$ in the infeasibility experiments.

\textbf{Time limit.} 
For the optimality experiments, we apply a time limit of 3600 seconds to small and medium-sized problems and 18000 seconds to large-scale problems for MIPLIB relaxation. A 15000-second time limit is used for Mittelmann’s benchmark.
For both infeasibility benchmark sets, the time limit is 3600 seconds per instance.

{
\textbf{Solver parameters.}
The termination and restart frequency are set to 200 iterations. For the primal weight update, we set $K_P=0.99$, $K_I=0.01$, $K_D=0.0$, and $\rho=0.3$ for the PID controller. 
For preconditioning, we sequentially perform 12 passes of geometric mean scaling, 10 iterations of Ruiz scaling, one $\ell_1$ scaling, and a final vector normalization. For restart, we set $\beta_{\mathrm{sufficient}}=0.2$, $\beta_{\mathrm{necessary}}=0.5$, $\beta_{\mathrm{artificial}}=0.36$.
For the active-set boost, we use the default settings: $W=10^4$, $\epsilon_r=\epsilon_c=10^{-8}$, $\epsilon_{\mathrm{act}}=10^{-4}$, $\kappa=1.1$, $\delta=0.05$, and $\theta=0.7$. The spectral norm is re-estimated when accumulated index additions and removals in $B$ and $F$ reach $1\%$ of the current $|B|+|F|$, and boosting is disabled after two rollbacks.
}

\textbf{Shifted geometric mean.} 
We use the shifted geometric mean of solving time to evaluate solver performance across a collection of instances. Formally, the shifted geometric mean is defined as $\left(\prod_{i=1}^n (t_i+\Delta)\right)^{1/n}-\Delta$, where $t_i$ is the solve time for the $i$-th instance. We shift by $\Delta = 10$ and denote this metric as SGM10. For the optimality experiments, an unsolved instance is assigned the corresponding time limit.

Tables~\ref{tab:miplib-no-presolve} and~\ref{tab:miplib-presolve} compare the performance of \texttt{cuPDLP.jl}, \texttt{HPR-lP}, \texttt{cuPDLP-C}, \texttt{Gurobi-PDLP} and \cupdlpx on 383 MIPLIB instances, with and without Gurobi presolve, respectively. Several key observations can be made:

\begin{itemize}
\item Across all four settings, \cupdlpx consistently outperforms the other solvers in terms of overall numerical performance, solving more instances and achieving notable reductions in running time.
\item The performance gain of \cupdlpx is especially significant for Small and Medium instances, while the improvement is less significant for Large instances, likely due to the high variance among the relatively small number of Large instances.
\end{itemize}

Table~\ref{tab:speedup} reports the speedup of \cupdlpx over \texttt{cuPDLP.jl}, \texttt{HPR-LP}, \texttt{cuPDLP-C}, and \texttt{Gurobi-PDLP}. For each baseline, the ``Hard'' rows consider instances on which either the baseline solver or \cupdlpx requires more than 10 seconds. Across all instances, \cupdlpx achieves speedups ranging from 1.3$\times$ to 3.6$\times$, while the speedups on hard instances range from 1.4$\times$ to 5.9$\times$. The largest improvements are observed over \texttt{cuPDLP.jl}, with speedups of 4.5$\times$--5.9$\times$ on hard instances. Overall, the gains are generally more pronounced on harder instances; for example, in the high-accuracy setting with presolve, \cupdlpx achieves a 3.0$\times$ speedup on hard instances and a 1.9$\times$ speedup overall relative to \texttt{cuPDLP-C}.

\begin{table}[h!]
\centering
\begin{tabular}{cccrcrcrcr}
\toprule
\multirow{2}{*}{$\epsilon$}     &    \multirow{2}{*}{\textbf{Solver}}                          & \multicolumn{2}{c}{\begin{tabular}[c]{@{}c@{}}\textbf{Small (269)} \\ (1-hour limit)\end{tabular}} & \multicolumn{2}{c}{\begin{tabular}[c]{@{}c@{}}\textbf{Medium (94)}\\ (1-hour limit)\end{tabular}}    & \multicolumn{2}{c}{\begin{tabular}[c]{@{}c@{}}\textbf{\textbf{Large (20)}}\\ (5-hour limit)\end{tabular}} & \multicolumn{2}{c}{\textbf{Total (383)}}          \\
                                                    & &\textbf{Count} & \textbf{Time}  & \textbf{Count} & \textbf{Time}  & \textbf{Count} & \textbf{Time} & \textbf{Count} & \textbf{Time}  \\ \midrule
\multirow{4}{*}{$\mathbf{10^{-4}}$} 
&\multicolumn{1}{c}{\textbf{cuPDLP.jl}}     & 266                   & 8.81                & 91                    & 11.55                & {19}                    & 77.43 &376 &11.07           \\
&\multicolumn{1}{c}{\textbf{HPR-LP}}          & 268 & 3.94  & 92 & 6.46   & 15 & 150.85 & 375 & 6.50 \\
&\multicolumn{1}{c}{\textbf{cuPDLP-C}}           & 268 & 3.78  & 92 & 7.54   & 18 & 60.12  & 378 & 5.92 \\

&\multicolumn{1}{c}{\textbf{Gurobi-PDLP}}        & 267 & 4.74  & 92 & 7.67   & 17 & 103.47 & 376 & 7.14 \\
&\multicolumn{1}{c}{\textbf{cuPDLPx} v0.1.0}     &  269                 &  2.77              &      94               &      4.62       &    19                & 57.18 & 382 & 4.39\\
\midrule
\multirow{4}{*}{$\mathbf{10^{-8}}$}
&\multicolumn{1}{c}{\textbf{cuPDLP.jl}}     & 260                   & 24.80               & 87                    & 36.90                & {17}                    & {208.66}  &364 &31.22             \\
&\multicolumn{1}{c}{\textbf{HPR-LP}}          & 261 & 10.77 & 91 & 16.78  & 14 & 485.71 & 366 & 16.09 \\
&\multicolumn{1}{c}{\textbf{cuPDLP-C}}           & 265 & 11.72 & 87 & 20.81  & 16 & 186.44 & 368 & 16.55 \\
&\multicolumn{1}{c}{\textbf{Gurobi-PDLP}}        & 261 & 12.22 & 89 & 19.77  & 17 & 388.87 & 367 & 17.76 \\
&\multicolumn{1}{c}{\textbf{cuPDLPx} v0.1.0}     &  266                 &   6.84             &     92                & 13.91            &  16                  &194.40  &374  &10.91 \\
\bottomrule
\end{tabular}
\caption{{SGM10 of solve time (in seconds)} on instances of {MIP Relaxations} without presolve.\label{tab:miplib-no-presolve}}
\end{table}

\begin{table}[h!]
\centering
\begin{tabular}{cccrcrcrcr}
\toprule
\multirow{2}{*}{$\epsilon$}     &   \multirow{2}{*}{\textbf{Solver}}                          & \multicolumn{2}{c}{\begin{tabular}[c]{@{}c@{}}\textbf{Small (269)} \\ (1-hour limit)\end{tabular}} & \multicolumn{2}{c}{\begin{tabular}[c]{@{}c@{}}\textbf{Medium (94)}\\ (1-hour limit)\end{tabular}}    & \multicolumn{2}{c}{\begin{tabular}[c]{@{}c@{}}\textbf{\textbf{Large (20)}}\\ (5-hour limit)\end{tabular}} & \multicolumn{2}{c}{\textbf{Total (383)}}          \\
                                                    & &\textbf{Count} & \textbf{Time}  & \textbf{Count} & \textbf{Time}  & \textbf{Count} & \textbf{Time} & \textbf{Count} & \textbf{Time}  \\ \midrule
\multirow{3}{*}{$\mathbf{10^{-4}}$}  
&\multicolumn{1}{c}{\textbf{cuPDLP.jl}}     & 268                   & 5.11                & 92                    & 9.03                & {19}                    & 26.74 &379 &6.75          \\
&\multicolumn{1}{c}{\textbf{HPR-LP}}          & 268 & 2.21  & 94 & 3.88   & 19 & 15.32  & 381 & 3.08 \\
&\multicolumn{1}{c}{\textbf{cuPDLP-C}}           & 268 & 1.99  & 94 & 4.25   & 19 & 13.25  & 381 & 2.95 \\
&\multicolumn{1}{c}{\textbf{Gurobi-PDLP}}        & 269 & 2.57  & 93 & 5.09   & 19 & 16.52  & 381 & 3.67 \\
&\multicolumn{1}{c}{\textbf{cuPDLPx} v0.1.0}& 269&	1.41&	94&	2.75&	19&	14.07&	382&	2.19\\
\midrule
\multirow{3}{*}{$\mathbf{10^{-8}}$}  
&\multicolumn{1}{c}{\textbf{cuPDLP.jl}}     & 264                   & 18.50             & 90                    & 29.40                & {19}                    & {63.68}  &373 &22.42             \\
&\multicolumn{1}{c}{\textbf{HPR-LP}}          & 268 & 5.85  & 92 & 11.90  & 19 & 45.63  & 379 & 8.32 \\
&\multicolumn{1}{c}{\textbf{cuPDLP-C}}           & 264 & 9.28  & 89 & 16.05  & 18 & 46.70  & 371 & 11.96 \\
&\multicolumn{1}{c}{\textbf{Gurobi-PDLP}}        & 265 & 7.35  & 90 & 14.52  & 19 & 39.65  & 374 & 9.95 \\
&\multicolumn{1}{c}{\textbf{cuPDLPx} v0.1.0}& 269&	3.98&	92&	9.60&	19&	43.09&	380&	6.28\\
\bottomrule
\end{tabular}
\caption{{SGM10 of solve time (in seconds)} on instances of {MIP Relaxations} with Gurobi presolve.\label{tab:miplib-presolve}}
\end{table}

\begin{table}[ht!]
\centering
\begin{tabular}{cccccc}
\toprule
\multirow{2}{*}{\textbf{Instance}} & \multirow{2}{*}{\textbf{Baseline}} 
& \multicolumn{2}{c}{\textbf{Without Presolve}} & \multicolumn{2}{c}{\textbf{With Presolve}} \\
\cmidrule(r){3-4} \cmidrule(l){5-6}
& & $10^{-4}$ & $10^{-8}$ & $10^{-4}$ & $10^{-8}$ \\
\midrule
\multirow{4}{*}{\textbf{Hard}}
  & cuPDLP.jl & 4.7 & 4.5 & 5.0 & 5.9 \\
  & HPR-LP & 2.1 & 1.9 & 1.6 & 1.5 \\
  & cuPDLP-C & 1.7 & 2.0 & 1.6 & 3.0 \\
  & Gurobi-PDLP     & 1.8 & 2.1 & 1.4 & 1.9 \\
\midrule
\multirow{4}{*}{\textbf{Overall}} 
  & cuPDLP.jl & 2.5 & 2.9 & 3.1 & 3.6 \\
  & HPR-LP & 1.5 & 1.5 & 1.4 & 1.3 \\
  & cuPDLP-C & 1.3 & 1.5 & 1.3 & 1.9 \\
  & Gurobi-PDLP     & 1.6 & 1.6 & 1.7 & 1.6 \\
\bottomrule
\end{tabular}
\caption{Summarization of \texttt{cuPDLPx} v0.1.0 speedup over \texttt{cuPDLP.jl}, \texttt{HPR-LP} \texttt{cuPDLP-C} and \texttt{Gurobi} on instances of MIP Relaxations. \label{tab:speedup}}
\vspace{4pt}
\begin{minipage}{\linewidth}\footnotesize\emph{Note.} The ``Hard'' rows report speedup on instances where either solver takes more than 10 seconds, while the ``Overall'' rows summarize speedup across all instances.\end{minipage}
\end{table}

Tables~\ref{tab:mittelmann-merged} and~\ref{tab:speedup-mittelmann} summarize the results on Mittelmann's LP benchmark set with 49 instances. Consistent with the MIPLIB experiments, \cupdlpx achieves consistent speedups over all other implementations under both tolerance settings, with and without presolve. Overall, the speedups range from 1.3$\times$ to 6.9$\times$ and are generally more pronounced on hard instances, with speedups of up to 16.1$\times$. For example, \cupdlpx achieves speedups of 1.3$\times$--1.7$\times$ over \texttt{HPR-LP} and 1.9$\times$--3.1$\times$ over \texttt{cuPDLP-C} across all instances, with the corresponding gains reaching 2.4$\times$ and 5.0$\times$, respectively, on hard instances.

\begin{table}[ht!]
\centering
\begin{tabular}{ccrcrcrcr}
\toprule
\multirow{3}{*}{\textbf{Solver}} 
& \multicolumn{4}{c}{\textbf{Without Presolve}} 
& \multicolumn{4}{c}{\textbf{With Presolve}} \\
\cmidrule(lr){2-5}\cmidrule(lr){6-9}
& \multicolumn{2}{c}{$\mathbf{10^{-4}}$} 
& \multicolumn{2}{c}{$\mathbf{10^{-8}}$} 
& \multicolumn{2}{c}{$\mathbf{10^{-4}}$} 
& \multicolumn{2}{c}{$\mathbf{10^{-8}}$} \\
& \textbf{Count} & \textbf{Time} 
& \textbf{Count} & \textbf{Time} 
& \textbf{Count} & \textbf{Time} 
& \textbf{Count} & \textbf{Time} \\
\midrule
\textbf{cuPDLP.jl}     & 43                     &   67.38             &  40                   &    219.51  & 46 & 48.76 & 42 & 187.22          \\
\textbf{HPR-LP} & 48 & 24.28 & 42 & 96.08 & 49 & 11.85 & 47 & 53.17 \\
\textbf{cuPDLP-C}    & 47 & 35.36 & 40 & 141.02 & 47 & 19.32 & 45 & 103.82 \\
\textbf{Gurobi-PDLP} & 44 & 32.47 & 39 & 153.94 & 46 & 21.95 & 43 & 93.29 \\
\textbf{cuPDLPx} v0.1.0     & 48 & 19.11 & 44 & 70.64  & 49 & 7.09  & 47 & 33.70 \\
\bottomrule
\end{tabular}
\caption{{SGM10 of solve time (in seconds)} on instances of Mittelmann's LP benchmark set with and without Gurobi presolve.\label{tab:mittelmann-merged}}
\end{table}

\begin{table}[ht!]
\centering
\begin{tabular}{cccccc}
\toprule
\multirow{2}{*}{\textbf{Instance}} & \multirow{2}{*}{\textbf{Baseline}} 
& \multicolumn{2}{c}{\textbf{Without Presolve}} & \multicolumn{2}{c}{\textbf{With Presolve}} \\
\cmidrule(r){3-4} \cmidrule(l){5-6}
& & $10^{-4}$ & $10^{-8}$ & $10^{-4}$ & $10^{-8}$ \\
\midrule
\multirow{4}{*}{\textbf{Hard}} 
  & cuPDLP.jl & 7.1 & 3.8 & 16.1 & 7.2 \\
  & HPR-LP & 1.3 & 1.5 & 2.4 & 1.8 \\
  & cuPDLP-C & 2.9 & 2.5 & 5.0 & 4.1 \\
  & Gurobi-PDLP     & 2.1 & 2.9 & 6.4 & 3.9 \\
\midrule
\multirow{4}{*}{\textbf{Overall}} 
  & cuPDLP.jl & 3.5 & 3.1 & 6.9 & 5.6 \\
  & HPR-LP & 1.3 & 1.4 & 1.7 & 1.6 \\
  & cuPDLP-C & 1.9 & 2.0 & 2.7 & 3.1 \\
  & Gurobi-PDLP     & 1.7 & 2.2 & 3.1 & 2.8 \\
\bottomrule
\end{tabular}
\caption{Summarization of \texttt{cuPDLPx} v0.1.0 speedup over \texttt{cuPDLP-C} and \texttt{Gurobi} on instances of Mittelmann's LP benchmark set. \label{tab:speedup-mittelmann}}
\vspace{4pt}
\begin{minipage}{\linewidth}\footnotesize\emph{Note.} The ``Hard'' rows report speedup on instances where either solver takes more than 10 seconds, 
while the ``Overall'' rows summarize speedup across all instances.\end{minipage}
\end{table}

{
\subsection{Comparison across different \cupdlpx versions}
\label{sec:version-comparison}
We compare five successive releases of \cupdlpx to evaluate the cumulative effects of algorithmic and implementation updates, using \texttt{v0.1.0} as the baseline.
\texttt{v0.2.3} integrates the PSLP presolver, \texttt{v0.2.5} further enables CUDA Graphs to reduce kernel-launch overhead, and \texttt{v0.3.0} adds support for \texttt{cuSPARSE} \texttt{cusparseSpMVOp} as an alternative SpMV backend.
Building on these updates, \texttt{v0.4.0} introduces active-set stepsize boosting, geometric mean scaling for diagonal preconditioning, and infeasibility detection. Earlier releases omit geometric mean scaling from the preconditioning sequence described above.

Tables~\ref{tab: cupdlpx 0.2.6 miplib} and~\ref{tab: cupdlpx 0.2.6 mittelmann-merged} compare these releases on the optimality benchmarks.
For the comparisons in this subsection, the reported solve time for \texttt{v0.2.3} and later versions includes the presolve time.
Experiments for \texttt{v0.3.0} are run with CUDA~13.3, which is required for \texttt{cusparseSpMVOp} support.
}

\begin{table}[h!]
\centering
{
\begin{tabular}{cccrcrcrcr}
\toprule
\multirow{2}{*}{}     &    \multirow{2}{*}{\textbf{\shortstack{cuPDLPx\\version}}}                          & \multicolumn{2}{c}{\begin{tabular}[c]{@{}c@{}}\textbf{Small (269)} \\ (1-hour limit)\end{tabular}} & \multicolumn{2}{c}{\begin{tabular}[c]{@{}c@{}}\textbf{Medium (94)}\\ (1-hour limit)\end{tabular}}    & \multicolumn{2}{c}{\begin{tabular}[c]{@{}c@{}}\textbf{\textbf{Large (20)}}\\ (5-hour limit)\end{tabular}} & \multicolumn{2}{c}{\textbf{Total (383)}}          \\
                                                    & &\textbf{Count} & \textbf{Time}  & \textbf{Count} & \textbf{Time}  & \textbf{Count} & \textbf{Time} & \textbf{Count} & \textbf{Time}  \\ \midrule
\multirow{5}{*}{$\mathbf{10^{-4}}$}
&\multicolumn{1}{c}{\textbf{v0.1.0}}     &  269                 &  2.77              &      94               &      4.62       &    19                & 57.18 & 382 & 4.39\\
&\multicolumn{1}{c}{\textbf{v0.2.3}}           & 269 & 1.67 & 94 & 3.52 & 19 & 19.42 & 382 & 2.69 \\
&\multicolumn{1}{c}{\textbf{v0.2.5}}           & 269 & 1.31 & 94 & 3.22 & 19 & 19.48 & 382 & 2.35 \\
&\multicolumn{1}{c}{\textbf{v0.3.0}}           & 269 & 1.05 & 94 & 2.96 & 19 & 16.95 & 382 & 2.04 \\
&\multicolumn{1}{c}{\textbf{v0.4.0}}           & 269 & 0.53 & 94 & 2.6 & 19 & 16.92 & 382 & 1.56 \\
\midrule
\multirow{5}{*}{$\mathbf{10^{-8}}$}
&\multicolumn{1}{c}{\textbf{v0.1.0}}     &  266                 &   6.84             &     92                & 13.91            &  16                  &194.40  & 374  & 10.91 \\
&\multicolumn{1}{c}{\textbf{v0.2.3}}           & 269 & 4.55 & 92 & 11.57 & 19 & 50.96 & 380 & 7.27 \\
&\multicolumn{1}{c}{\textbf{v0.2.5}}           & 269 & 3.65 & 93 & 10.99 & 19 & 50.77 & 381 & 6.4 \\
&\multicolumn{1}{c}{\textbf{v0.3.0}}           & 269 & 2.89 & 93 & 9.65 & 19 & 45.65 & 381 & 5.43 \\
&\multicolumn{1}{c}{\textbf{v0.4.0}}           & 269 & 1.64 & 93 & 7.09 & 19 & 34.09 & 381 & 3.71 \\
\bottomrule
\end{tabular}
}
\caption{SGM10 of solve time (in seconds) on instances of {MIP Relaxations} for different \cupdlpx versions.\label{tab: cupdlpx 0.2.6 miplib}}
\end{table}

\begin{table}[ht!]
\centering
{
\begin{tabular}{ccrcrcrcr}
\toprule
\multirow{2}{*}{\textbf{\shortstack{cuPDLPx\\version}}}
& \multicolumn{2}{c}{$\mathbf{\epsilon=10^{-4}}$} 
& \multicolumn{2}{c}{$\mathbf{\epsilon=10^{-8}}$} \\
& \textbf{Count} & \textbf{Time} 
& \textbf{Count} & \textbf{Time} \\
\midrule
\textbf{v0.1.0}     & 48 & 19.11 & 44 & 70.64   \\
\textbf{v0.2.3}     & 48 & 12.06 & 47 & 43.06   \\
\textbf{v0.2.5}     & 48 & 11.03 & 47 &  39.36  \\
\textbf{v0.3.0}     & 49 & 9.62 & 47 &  34.29  \\
\textbf{v0.4.0}     & 49 & 6.84 & 49 &  19.03  \\
\bottomrule
\end{tabular}
}
\caption{SGM10 of solve time (in seconds) on instances of Mittelmann's LP benchmark set for different \cupdlpx versions.\label{tab: cupdlpx 0.2.6 mittelmann-merged}}
\end{table}

Across both benchmark sets and tolerance settings, successive updates give \texttt{v0.4.0} a $2.8\times$--$3.7\times$ speedup over \texttt{v0.1.0}, while solving more instances at the tighter tolerance.
Specifically, the introduction of PSLP in \texttt{v0.2.3} reduces SGM10 across all three MIP size groups, with the largest relative gains on large instances.
The further gains from \texttt{v0.4.0} over \texttt{v0.3.0} are more pronounced at $10^{-8}$, consistent with the active-set stepsize boost being applied only in these runs due to its $10^{-4}$ activation tolerance.
In addition, \cupdlpx \texttt{v0.4.0} incurred no false infeasibility detections in these optimality experiments.

Relative to \texttt{cuPDLP.jl} without external presolve (Tables~\ref{tab:miplib-no-presolve} and~\ref{tab:mittelmann-merged}), \cupdlpx \texttt{v0.4.0} achieves speedups of $7.1\times$--$8.4\times$ on MIP Relaxations and $9.9\times$--$11.5\times$ on Mittelmann's LP benchmark across the two tolerance settings (with its integrated PSLP presolve time included).

\subsection{Infeasibility detection}
\label{sec:infeasibility-experiments}
Tables~\ref{tab: miplib infeas} and~\ref{tab: netlib infeas} report the results on MIPLIB 2017 Infeas and the infeasible Netlib instances, respectively, with and without PSLP presolve~\cite{Cederberg26}.

The tables report the number of accepted infeasibility or unboundedness detections (Count) and SGM10 in seconds (Time).
A detection is counted only if the reported status is consistent with the known infeasibility or unboundedness of the instance. Runs that terminate with an optimality status are excluded from Count.

\begin{table}[ht!]
\centering
\begin{tabular}{ccrcr}
\toprule
\multirow{2}{*}{\textbf{Solver}}
& \multicolumn{2}{c}{\textbf{Without Presolve}}
& \multicolumn{2}{c}{\textbf{With PSLP}} \\
\cmidrule(lr){2-3}\cmidrule(lr){4-5}
& \textbf{Count} & \textbf{Time} & \textbf{Count} & \textbf{Time} \\
\midrule
\textbf{Gurobi} primal simplex      & 323 & 11.52  & 334 & 4.55   \\
\textbf{Gurobi} dual simplex        & 332 & 6.05   & 332 & 2.98   \\
\textbf{Gurobi} barrier             & 187 & 153.63 & 253 & 39.65  \\
\textbf{Gurobi} homogeneous barrier & 257 & 42.84  & 288 & 15.59  \\
\textbf{HiGHS}                      & 330 & 14.47  & 330 & 6.23   \\
\textbf{SCS}                        & 187 & 264.62 & 246 & 84.00  \\
\textbf{OR-Tools PDLP}              & 243 & 305.34 & 277 & 103.99 \\
\textbf{H-PDLP}                    & 326 & 6.17   & 334 & 2.10   \\
\textbf{cuPDLPx} v0.4.0            & 333 & 3.45   & 334 & 1.31   \\
\bottomrule
\end{tabular}
\caption{Infeasibility detection results on the 334 MIPLIB 2017 Infeas instances, with and without PSLP presolve.\label{tab: miplib infeas}}
\end{table}

\begin{table}[ht!]
\centering
\begin{tabular}{ccrcr}
\toprule
\multirow{2}{*}{\textbf{Solver}}
& \multicolumn{2}{c}{\textbf{Without Presolve}}
& \multicolumn{2}{c}{\textbf{With PSLP}} \\
\cmidrule(lr){2-3}\cmidrule(lr){4-5}
& \textbf{Count} & \textbf{Time} & \textbf{Count} & \textbf{Time} \\
\midrule
\textbf{Gurobi} primal simplex      & 29 & 0.01  & 29 & 0.01  \\
\textbf{Gurobi} dual simplex        & 29 & 0.04  & 29 & 0.00  \\
\textbf{Gurobi} barrier             & 23 & 24.01 & 27 & 2.28  \\
\textbf{Gurobi} homogeneous barrier & 27 & 5.05  & 29 & 0.01  \\
\textbf{HiGHS}                      & 27 & 5.07  & 29 & 0.00  \\
\textbf{SCS}                        & 20 & 3.93  & 26 & 0.61  \\
\textbf{OR-Tools PDLP}              & 24 & 74.55 & 27 & 21.97 \\
\textbf{H-PDLP}                    & 29 & 2.01  & 29 & 1.21  \\
\textbf{cuPDLPx} v0.4.0            & 29 & 1.52  & 29 & 0.88  \\
\bottomrule
\end{tabular}
\caption{Infeasibility detection results on the 29 infeasible Netlib instances, with and without PSLP presolve.\label{tab: netlib infeas}}
\end{table}

On MIPLIB 2017 Infeas, \cupdlpx is the fastest solver overall in both presolve settings. It detects infeasibility or unboundedness on 333 of 334 instances without presolve and all 334 instances with PSLP.
The remaining instance without presolve terminates with an optimality status and is not counted as a successful infeasibility detection.
On Netlib, \cupdlpx detects infeasibility on all 29 instances in both settings.
Gurobi's primal and dual simplex methods are faster overall on this benchmark and also detect infeasibility on every instance.
This is consistent with the small size of the Netlib instances, which are relatively easy for simplex methods.

\section{Conclusions}
In this paper, we present \cupdlpx, an enhanced GPU implementation in C and CUDA of restarted Halpern PDHG for solving large-scale LPs. The solver combines fixed-point-error-based restart criteria, a PID-controlled primal weight update, active-set stepsize boosting, and infeasibility detection with a highly optimized GPU implementation and integration of the PSLP presolver. Compared with \texttt{cuPDLP.jl}, \cupdlpx \texttt{v0.4.0} achieves speedups of $7.1\times$--$8.4\times$ on MIPLIB 2017 LP relaxations and $9.9\times$--$11.5\times$ on Mittelmann's LP benchmark across different tolerances. On MIPLIB 2017 Infeas, it detects infeasibility or unboundedness on all 334 instances with PSLP while incurring no false infeasibility detections on any feasible benchmark instance. These results highlight the effectiveness of algorithmic and implementation enhancements in exploiting GPUs for high-performance linear programming.

\section*{Acknowledgement}
Haihao Lu is partially supported by AFOSR Grant No. FA9550-24-1-0051 and ONR Grant No. N000142412735. Zedong Peng is partially supported by ONR Grant No. N000142412735. Jinwen Yang is partially supported by AFOSR Grant No. FA9550-24-1-0051. The authors gratefully acknowledge the support of the NVIDIA Academic Grant Program.

\bibliographystyle{amsplain}
\bibliography{ref-papers-v5.0}

\end{document}